\documentclass[12pt,a4paper,twoside]{article}
\usepackage{graphicx}
\usepackage{times}
\usepackage{mathptmx}
\usepackage{caption,subfig}
\usepackage{tikz}
\usepackage{caption}

\usepackage{cite}
\usepackage[T1,OT1]{fontenc}
\usepackage{textcomp}
\usepackage{xcolor}
\usepackage{graphicx}

\usepackage[colorlinks,
linkcolor=blue,
anchorcolor=blue,
citecolor=blue,]{hyperref}

\usepackage{multirow}
\usepackage{mathrsfs,amssymb,amsthm,stmaryrd,amsmath,latexsym,indentfirst}
\usepackage{stmaryrd}
\usepackage{makecell}
\usepackage{booktabs}
\usepackage{xcolor}
\usepackage{subfig}
\usepackage{bm}
\usepackage[ruled,linesnumbered,vlined]{algorithm2e}
\usepackage[symbol]{footmisc}
\usepackage{cellspace}

\usepackage[capitalise]{cleveref}
\usepackage{chngcntr}

\usepackage{amsmath}
\numberwithin{equation}{subsection}

\newtheorem{mytheorem}{Theorem}[section]
\newtheorem{mylemma}[mytheorem]{Lemma}
\newtheorem{mycorollary}[mytheorem]{Corollary}

\newtheorem{myconj}[mytheorem]{Conjecture}
\newtheorem{mycase}{Case}[mytheorem]

\newtheorem{submycase}{Case}[mycase]

\theoremstyle{definition}
\newtheorem{myremark}{Remark}

\newtheorem{claim}{Claim}[subsection]
\renewenvironment{proof}[1][Proof]{%
	\par\noindent{\normalfont\bfseries #1.}\enspace%
	\setcounter{mycase}{0}
	\setcounter{submycase}{0}
	\ignorespaces
}{%
	\unskip\hfill$\square$\par
}
\newenvironment{claimproof}
{\par\noindent \textit{Proof.}}
{\qed\par}
\usepackage{hyperref} 

\begin{document}
	
	\title{{Spanning cyclability of split graphs}
		}
	\author{Jingfan Zai, Wenjing Lv, Eminjan Sabir \footnote{This research was funded by Natural Science Foundation of Xinjiang Uygur Autonomous Region (No. 2025D01C34) and National Natural Science Foundation of China (No. 12261085).}\thanks{Corresponding author: eminjan20150513@163.com}}
	\date{ College of Mathematics and System Sciences, Xinjiang University, \\Urumqi,
		830046, P. R. China}
	\maketitle
	\renewcommand{\abstractname}{}
	\begin{abstract}
		\noindent
		{\bf Abstract:} {  For a positive integer $k$, a graph $G$ is called \textbf{$k$-spanning cyclable} (resp. \textbf{$k$-edge spanning cyclable}) if, given any $k$ distinct vertices $v_1,v_2,\dots,v_k$ (resp. $k$ independent edges $e_1,e_2,\dots,e_k$) of $G$, there exist $k$ vertex-disjoint cycles $C_1,C_2,\dots,C_k$ in $G$ such that $v_i\in V(C_i)$ (resp. $e_i\in E(C_i)$) for each $i$ and $V(C_1)\cup V(C_2)\cup\cdots\cup V(C_k)=V(G)$. In particular, the notion of $1$-spanning cyclable coincides with the existence of a Hamiltonian cycle, and a Hamilton-connected graph is necessarily $1$-edge spanning cyclable. For a positive integer $ r\ge3$, a graph $G$ is \textbf{$K_{1,r}$-free} if $G$ does not have an induced subgraph isomorphic to $K_{1,r}$, and $K_{1,4}+e$ denotes the graph obtained by adding an edge $e$ between two nonadjacent vertices of $K_{1,4}$. There are several interesting conjectures and sparse results relative to these conjectures on Hamiltonian properties of $K_{1,r}$-free \textbf{split graphs}, whose vertex set can be partitioned into a maximum clique $W$ and an independent set $J$. In this paper, we extend some of these results by proving the following:

$(1)$ every $(k+1)$-connected $K_{1,3}$-free split graph of order at least $4k-1$ is $k$-spanning cyclable for $k\ge 1$; in particular, every $2$-connected $K_{1,3}$-free split graph of order at least $3$ is Hamiltonian;

$(2)$ every $2k$-connected $K_{1,3}$-free split graph is $k$-edge spanning cyclable for $k\ge 2$;

$(3)$ every $(k+1)$-connected $\{K_{1,4},K_{1,4}+e\}$-free split graph with $|W|\ge |J|\ge 3k+5$ is $k$-spanning cyclable for $k\ge 1$; in particular, every $2$-connected $\{K_{1,4},K_{1,4}+e\}$-free split graph with $|W|\ge |J|\ge 8$ is Hamiltonian;

$(4)$ every $2k$-connected $\{K_{1,4},K_{1,4}+e\}$-free split graph with order at least $4k+9$ is $k$-edge spanning cyclable for $k\ge 2$.}

		\begin{flushleft}
			\textbf{Keywords:}  Split graphs; $K_{1,r}$-free graphs; Spanning cyclable; Edge spanning cyclable; Hamiltonian graphs
		\end{flushleft}
		
	\end{abstract}
	

	\section{Introduction} \label{Sec.1}

	In this paper, we only consider finite, connected and simple graphs. Let $G$ be a graph with vertex set $V(G)$ and edge set $E(G)$. A graph $G$ is \textbf{Hamiltonian} if it contains a spanning cycle and $G$ is \textbf{Hamilton-connected} if for any two distinct vertices $u,v\in V(G)$, there exists a spanning $(u,v)$-path in $G$.  In addition, $\delta(G)$ and $\kappa(G)$ represent the minimum degree and connectivity of $G$. For a cycle $C$ in the graph $G$, we write $C = v_1v_2\cdots v_kv_1$, which means vertices $v_1,v_2,\dots,v_k$ are adjacent in order and form a cycle with $v_k$ connected back to $v_1$. For any positive integer $k$, let $[k] = \{1,2,\dots,k\}$ denote the set of all positive integers from $1$ to $k$. For any edge $e \in E(G)$, let $V(e)$ denote the set of its two vertices.

	A graph $G$ is called \textbf{$k$-spanning cyclable} if for any $k$ distinct vertices $v_1,v_2,\dots,v_k$ of $G$, there exist $k$ vertex-disjoint cycles $C_1,C_2,\dots,C_k$ in $G$ such that $v_i\in V(C_i)$ for each $i$ and $V(C_1)\cup V(C_2)\cup\cdots\cup V(C_k)=V(G)$\cite{13}. Similarly, a graph $G$ is called \textbf{$k$-edge spanning cyclable} if for any $k$ independent edges $e_1,e_2,\dots,e_k$ of $G$, there exist $k$ vertex-disjoint cycles $C_1,C_2,\dots,C_k$ in $G$ such that $e_i\in E(C_i)$ for each $i$ and $V(C_1)\cup V(C_2)\cup \cdots \cup V(C_k)=V(G)$ \cite{22}. The $k$-spanning cyclability and $k$-edge-spanning cyclability of graphs have close connections with Hamiltonian properties. In particular, the notion of $1$-spanning cyclable coincides with the existence of a Hamiltonian cycle, and a Hamilton-connected graph is necessarily $1$-edge spanning cyclable. Recent progress on the spanning cyclability of graphs can be found in \cite{11,12,14,23,24,25}. In particular, Yang et al.\cite{10} gave a necessary minimum degree condition for spanning cyclable graphs.
	
	\begin{mytheorem}(Yang et al.\cite{10})\label{Th1.1}
		If a graph $G$ is $k$-spanning cyclable, then every vertex of $G$ has degree at least $k+1$.
	\end{mytheorem}
	
For a positive integer $ r\ge3$, a graph $G$ is \textbf{$K_{1,r}$-free} if $G$ does not have an induced subgraph isomorphic to $K_{1,r}$. A connected graph $G$ is called a \textbf{split graph} if its vertex set can be partitioned into two disjoint subsets $W$ and $J$ (either of which may be empty), where $W$ is a  maximum clique and $J$ is an independent set. Split graphs were introduced by Foldes and Hammer \cite{02} in 1977 and have been further studied in \cite{03, 04, 05, 06, 07,19,20,26,27,28,29,30,31,32,33}. Recently, the necessary and sufficient conditions for Hamiltonian properties of $K_{1,r}$-free split graphs have been well-studied.
	
	\begin{mytheorem}(Renjith and Sadagopan \cite{05}) \label{Theorem 1.1}
		 Let $G$ be a $K_{1,3}$-free split graph. Then $G$ is Hamiltonian if and only if $G$ is $2$-connected.
		
	\end{mytheorem}
	
	Let $G$ be a $K_{1,4}$-free split graph where $W$ is a clique and $J$ is an independent set. Consider a subgraph $H$ of $G$ where $V_a=\{x\in J\mid d_G(x)=2\}$, $V_b=N_G(V_a)$, $V(H)=V_a\cup V_b$, and $E(H)=\{xy\in E(G)\mid x\in V_a,y\in V_b\}$. Then $H$ is a bipartite subgraph of $G$. Let $C$ be an induced cycle of $H$ satisfying $W\setminus V(C)\neq\emptyset$. We call such a cycle $C$ \textit{a short cycle} of $G$.
	
	\begin{mytheorem}(Renjith and Sadagopan \cite{05}) \label{Theorem 1.2}
		Let $G$ be a $2$-connected $K_{1,4}$-free split graph with \(|W|\ge |J|\ge 8\). Then $G$ is Hamiltonian if and only if there is no short cycle in $G$.
		
	\end{mytheorem}
	
	In 2022, Dai et al. \cite{07} further extended the above results and proposed a conjecture that $K_{1,r}$-free split graphs are Hamiltonian.
	
	\begin{mytheorem}(Dai et al. \cite{07}) \label{Theorem 1.3}
		Let $G$ be a $K_{1,3}$-free split graph. Then $G$ is pancyclic if and only if $G$ is $2$-connected.
		
	\end{mytheorem}
	
	\begin{mytheorem}(Dai et al. \cite{07}) \label{Theorem1.4}
		Let $G$ be a $K_{1,4}$-free split graph. Then $G$ is Hamiltonian if and only if $G$ is $3$-connected.
		
	\end{mytheorem}

	\begin{myconj}(Dai et al. \cite{07}) \label{Conjecture 1.5}
		 Let $G$ be a $K_{1,r+1}$-free split graph with at least three vertices. If $G$ is $r$-connected, then $G$ is Hamiltonian.
		
	\end{myconj}
	
	In 2023, Liu et al. \cite{06} further investigated the Hamilton-connectedness of $K_{1,r}$-free split graphs and proposed a conjecture.
	
	\begin{mytheorem}(Liu et al. \cite{06}) \label{Theorem1.6}
		 For any $ r \in \{3,4\} $,  if $G$ is an $r$-connected $K_{1,r}$-free split graph, then $G$ is Hamilton-connected.
	\end{mytheorem}
	
	\begin{myconj}(Liu et al. \cite{06}) \label{Conjecture 1.7}
		 Let $r\ge 3$ be an integer. Every $r$-connected $K_{1,r}$-free split graph is Hamilton-connected.
		
	\end{myconj}
	
	 Let $r\ge 3$ be an integer, and let $K_{1,r+1}+e$ denote the graph obtained by adding an edge between two nonadjacent vertices of $K_{1,r+1}$. Obviously, a $K_{1,3}$-free graph is also a $\{K_{1,4}, K_{1,4}+e\}$-free graph, but the reverse is not true. The Hamiltonicity of $\{K_{1,4},K_{1,4}+e\}$-free graphs was first investigated by Li and Schelp \cite{08} in $2002$, and Li \cite{09} studied the Hamiltonicity of $2$-connected graphs in 2004.
	
	\begin{mytheorem}(Li \cite{09})\label{Theorem1.8}
		Let $G$ be a $2$-connected $\{K_{1,4},K_{1,4}+e\}$-free graph with at least $13$ vertices. If $\delta(G)\ge \frac{n}{4}$, then $G$ is Hamiltonian or $G\in \mathcal{F}$, where $\mathcal{F}$ is a family of non-Hamiltonian graphs of connectivity $2$.
	\end{mytheorem}
	
	Ryj\'a\v{c}ek et al. \cite{18} proved that the minimum degree condition can be significantly decreased to $6$ if the connectivity is increased to at least $5$ when considering the Hamiltonicity of $\{K_{1,4}, K_{1,4}+e\}$-free graphs. It is well known that the line graph does not contain $K_{1,3}$ as an induced subgraph. Thus Ryj\'a\v{c}ek et al. \cite{18} also posed the following two equivalent conjectures.
	
	\begin{mytheorem}(Ryj\'a\v{c}ek et al. \cite{18}) \label{Theorem1.9}
		Every $5$-connected $\{K_{1,4}, K_{1,4}+e\}$-free graph with minimum degree at least $6$ is Hamiltonian.
	\end{mytheorem}
	
	\begin{myconj}(Ryj\'a\v{c}ek et al. \cite{18}) \label{Conjecture 1.10}
		Every $4$-connected $\{K_{1,4}, K_{1,4}+e\}$-free graph with minimum degree at least $6$ is Hamiltonian.
		
	\end{myconj}
	
	\begin{myconj}(Ryj\'a\v{c}ek et al. \cite{18}) \label{Conjecture 1.11}
		Every $4$-connected line graph with minimum degree at least $5$ is Hamiltonian.
	\end{myconj}
	
	There exists an infinite family of $r$-connected $K_{1,r+1}$-free split graphs, each of which is not Hamilton-connected \cite{06}. Thus, Tian et al. \cite{19} focused on Hamilton-connectedness of $\{K_{1,4}, K_{1,4}+e\}$-free split graphs, and proposed an interesting conjecture.
	
	\begin{mytheorem}(Tian et al. \cite{19}) \label{Theorem1.12}
		If $G$ is a $3$-connected $\{K_{1,4}, K_{1,4}+e\}$-free split graph with vertices at least $13$, then $G$ is Hamilton-connected.
	\end{mytheorem}
	
	\begin{myconj}(Tian et al. \cite{19}) \label{Conjecture 1.13}
		Let $r \ge 3$ be an integer. Every $r$-connected $\{K_{1,r+1}, K_{1,r+1}+e\}$-free split graph with vertices at least $2r+7$ is Hamilton-connected.
	\end{myconj}
	
	This paper mainly investigates the spanning cyclability and edge spanning cyclability of $K_{1,3}$-free and $\{K_{1,4},K_{1,4}+e\}$-free split graphs.

\begin{mytheorem}\label{Theorem1.14}
		Let $G$ be a \(K_{1,3}\)-free split graph with a split partition \((W,J)\).
	If \(|V(G)| \ge 4k-1\) and \(G\) is \((k+1)\)-connected for $k\ge 1$, then \(G\) is \(k\)-spanning cyclable.

\end{mytheorem}

	\begin{myremark}
	It follows from Theorem \ref{Th1.1} that the $(k+1)$-connectedness condition in Theorem \ref{Theorem1.14} is sharp.
\end{myremark}

\begin{mycorollary}
	Let $G$ be a $K_{1,3}$-free split graph with a split partition $(W,J)$. If $|V(G)|\ge 3$ and $G$ is $2$-connected, then $G$ is Hamiltonian.
\end{mycorollary}

\begin{mytheorem} \label{Theorem1.16}
	Let \(G\) be a \(K_{1,3}\)-free split graph with a split partition \((W,J)\).
	If \(G\) is \(2k\)-connected for $k \ge 2$, then \(G\) is \(k\)-edge spanning cyclable.
\end{mytheorem}

\begin{myremark}
	The $2k$-connectedness bound in Theorem \ref{Theorem1.16} is sharp.
	To see this, take a vertex $a \in J$ with $d_G(a) = 2k-1$ in a $(2k-1)$-connected $K_{1,3}$-free split graph. Let $e_k=ab$ with $b\in N_G(a)$, and take $k-1$ vertex-disjoint edges $e_1,\dots,e_{k-1}$ satisfying $V(e_i)\subseteq N_G(a)\setminus\{b\}$ for each $i\in[k-1]$.
	These $k-1$ edges occupy $2k-2$ distinct neighbors of $a$, leaving only $b$. No extra neighbor of $a$ is available to form a cycle $C_k$ containing $e_k$.
	Hence this graph is not $k$-edge spanning cyclable.
\end{myremark}

\begin{mytheorem} \label{Theorem1.17}
	Let \(G\) be a \(\{K_{1,4},K_{1,4}+e\}\)-free split graph with a split partition \((W,J)\).
	If  \(|W|\ge |J|\ge 3k+5\) and \(G\) is \((k+1)\)-connected for $ k \ge 1$, then \(G\) is \(k\)-spanning cyclable.
	
\end{mytheorem}

\begin{mycorollary}
	Let $G$ be a $\{K_{1,4}, K_{1,4}+e\}$-free split graph with a split partition $(W,J)$. If $|W|\ge |J|\ge 8$ and $G$ is $2$-connected, then $G$ is Hamiltonian.
\end{mycorollary}

	\begin{mytheorem}\label{Theorem1.19}
	Let \(G\) be a \(\{K_{1,4},K_{1,4}+e\}\)-free split graph with a split partition \((W,J)\).
	If  \(|V(G)|\ge 4k+9\) and \(G\) is \(2k\)-connected for $ k \ge 2$, then \(G\) is \(k\)-edge spanning cyclable.
\end{mytheorem}

 The remainder of this paper is organized as follows. Section \ref{Sec.2} summarizes the lemmas used in subsequent arguments, and Section \ref{Sec.3} provides the proofs of our main results.

	\section{Preliminaries} \label{Sec.2}
	
	In this paper, any undefined terms and notation can be found in \cite{01}. Throughout this section, we assume that $J \neq \emptyset$ and let $r \ge 3$ be an integer. For a vertex $v \in V(G)$, we define the neighborhood of $v$ in $G$ as
	\[
	N_G(v) = \{z \in V(G) : vz \in E(G)\}.
	\]
	Given a split graph $G$ with the split partition $(W,J)$, for a vertex $v \in W$, we define
	\[
	N_J(v) = N_G(v) \cap J,\quad \text{and } W_J = W \cap N_G(J).
	\]
	
	For any vertex $j \in J$, since $W$ is a maximum clique of $G$, we conclude that $N_G(j)$ is a proper subset of $W$, and hence is a vertex cut of $G$. This leads to the following.
	\[
	|W| \ge |N_G(j)| + 1 \ge \kappa(G) + 1. \tag{1} \label{eq1}
	\]

	Since $G$ is a split graph with $G[W]$ being a complete graph, if $J_0 \ne \emptyset$, then every vertex cut either equals $W$, or has the form $N_G(j)$ for a vertex $j \in J$. It follows from (\ref{eq1}) that
	\[
	\kappa(G) = \min\{|N_G(j)| : j \in J\}. \tag{2} \label{eq2}
	\]
	
	Since $(W,J)$ is the split partition of $G$ and by (\ref{eq2}), we have the following lemmas.
	
	\begin{mylemma}(Xiong et al. \cite{20}) \label{lem2.1}
		Let $G$ be a $K_{1,r}$-free split graph with the split partition $(W,J)$ and $\kappa(G) \ge 2$, and let $j \in J$ be a vertex. Each of the following holds.
		\begin{enumerate}
			\item[(i)] $G-j$ is also a $K_{1,r}$-free split graph with the split partition $(W, J - \{j\})$.
			\item[(ii)] $\kappa(G-j) \ge \kappa(G)$.
		\end{enumerate}
	\end{mylemma}
		
		\begin{mylemma}(Xiong et al. \cite{20}) \label{lem2.2}
			Let $G$ be a $K_{1,r}$-free split graph with $\kappa(G) \ge 2$ and with the split partition $(W,J)$ such that $J \ne \emptyset$. Let $J_0 \subseteq J$ be the vertex subset where each vertex $j_0 \in J_0$ if and only if $d_G(j_0) = |W| - 1$. For any $w \in W - W_J$, each of the following holds.
			\begin{enumerate}
				\item[(i)] If $J_0 = \emptyset$, then for every $j \in J$, $d_G(j) < |W| - 1$. In this case, $G - w$ is also a $K_{1,r}$-free split graph with the split partition $(W - \{w\}, J)$.
				\item[(ii)] If $J_0 \ne \emptyset$, then as $w \in W - W_J$, every vertex $w' \in W - \{w\}$ is adjacent to a vertex in $J_0$. In this case, $G - w$ is also a $K_{1,r}$-free split graph with the split partition $((W - \{w\}) \cup J_0, J - J_0)$.
				\item[(iii)] In any case, $|W| - 1 \ge \kappa(G - w) = \kappa(G)$.
			\end{enumerate}
		\end{mylemma}
		
		By Lemmas \ref{lem2.1} and \ref{lem2.2}, if $G$ is a $K_{1,r}$-free split graph with split partition $(W,J)$, then for any vertex $v \in V(G)$, $G-\{v\}$ is still a $K_{1,r}$-free split graph. Moreover, the connectivity of $G-\{v\}$ would be reduced only when $v \in W_J$.

	\section{Proof of Main Results } \label{Sec.3}

	\subsection{Spanning Cyclability of $K_{1,3}$-Free Split Graphs}

\vskip 0.3cm
	
	\begin{proof}[Proof of Theorem \ref{Theorem1.14}]
		Let G be a $(k+1)$-connected \(K_{1,3}\)-free split graph with a split partition \((W,J)\). The proof is by induction on k. The base case for $k=1$ holds by Theorem \ref{Theorem 1.1}.
		Now assume that for  $k \ge 2$, every $k$-connected graph $G$ with $|V(G)| \ge 4k-5$ is $(k-1)$-spanning cyclable. Let $v_1,v_2,\dots,v_k$ be any $k$ vertices of $G$ and set $u=v_k$.
	
	\begin{claim} \label{claim3}
		Assume $s \in J$, and there exist $x, y \in  N_G(s) $ such that for any $j \in J \setminus \{s\}$, if $|N_G(j) \cap \{x, y\}| = 2$,
		then $\max\{d_G(s), d_G(j)\}\ge k+2$.
	\end{claim}
	
	\begin{claimproof}
		We proceed by contradiction. Suppose $d_G(s) = d_G(j) = k+1$. Since $G$ is $K_{1,3}$-free and $W$ is a maximal clique, any vertex in $W$ that is non-adjacent to $s$ must be adjacent to $j$. Hence
		\begin{equation}\label{eq:1}
			|W| \le (k+1)+(k+1)-2 = 2k.
		\end{equation}
		Since $|V(G)| = |W| + |J|$, we have
		\begin{equation}\label{eq:2}
			|J| = |V(G)| - |W| \ge 4k-1 - 2k = 2k-1.
		\end{equation}
		Since $ (k+1)|J| \le e(W,J)\le 2|W|$, we obtain
		\begin{equation}\label{eq:3}
			(k+1)|J| \le 2|W|.
		\end{equation}
		Combining \eqref{eq:1}, \eqref{eq:2} and \eqref{eq:3}, we have
		\[
		(k+1)(2k-1) \le 4k.
		\]
		This inequality contradicts the assumption $k \ge 2$.
	\end{claimproof}
	
	\begin{mycase}\label{3.2case1}
		For any $t\in W$, $|N_J(t)|\le 1$.
	\end{mycase}
	
	\begin{submycase}\label{3.2case1.1}
		$u\in J$.
	\end{submycase}
	
	Since $G$ is $(k+1)$-connected, there must exist two vertices $x,y\in N_G(u)\setminus\{v_1,v_2,\dots,v_{k-1}\}$. Let $C_k=uxyu$ and set $G'=G-C_k$. Then $G'$ is a $k$-connected \(K_{1,3}\)-free split graph.

		\begin{submycase}\label{3.2case1.2}
		$u,v_1,v_2,\dots,v_{k-1}\in W$.
	\end{submycase}
	
	If $N_J(u)\neq\emptyset$, then there exists $z\in N_J(u)$. As $d_G(z)\ge k+1$, we may choose $a\in N_G(z)\setminus\{u,v_1,v_2,\dots,v_{k-1}\}$ and set $C_k=uazu$.
	If $N_J(u)=\emptyset$, the $(k+1)$-connectedness of $G$ and the maximality of $W$ give $|W| \ge k+2$. Then there must exist two vertices $x,y\in W\setminus\{u,v_1,v_2,\dots,v_{k-1}\}$.  If $N_J(x)\cap N_J(y)=\emptyset$, we take $C_k=uxyu$. If $N_J(x)\cap N_J(y)\neq\emptyset$, we may assume $N_J(x)\cap N_J(y)=\{j\}$ and take $C_k=uxjyu$.
	Set $G'=G-C_k$, then $G'$ is a $k$-connected \(K_{1,3}\)-free split graph.
	
	\begin{mycase}\label{3.2case2}
		There is some $t\in W$ with $|N_J(t)| \ge 2$.
	\end{mycase}
	
	As $G$ is $K_{1,3}$-free, every vertex in $W$ has at most $2$ neighbors in $J$, so $|N_J(t)| = 2$.
	Moreover, every vertex in $W$ is adjacent to at least one vertex in $J$, otherwise, a $K_{1,3}$ would be induced.
		
		\begin{submycase}\label{3.2case2.1}
		$u\in J$.
	\end{submycase}
	
	By the $(k+1)$-connectedness of $G$, there must exist two vertices $x,y\in N_G(u)\setminus\{v_1,v_2,\dots,v_{k-1}\}$. Excluding the vertex $u$, if $x$ and $y$ have no other common neighbors in $J$, we take $C_k=uxyu$. Otherwise, let $j$ be an additional common neighbor of $x$ and $y$ in $J$. Then we take $C_k=uxjyu$ provided that $j\neq v_i$ for any $1\le i\le k-1$.
	Now consider the case $j\in\{v_1,v_2,\dots,v_{k-1}\}$, and without loss of generality let $j=v_1$. In this case, by Claim \ref{claim3}, $\max\{d_G(u), d_G(v_1)\}\ge k+2$. Therefore, if $d_G(v_1)\ge k+2$, we may take $C_k=uxyu$.
		
	\begin{claim} \label{claim4}
		If $d_G(u) \ge k+2$ and $d_G(v_1) = k+1$, then there exists a vertex $z \notin \{ v_2, \dots, v_{k-1}, x, y\}$ such that $z \in N_G(u) $ and $z \notin N_G(v_1)$.
	\end{claim}
	
	\begin{claimproof}
		We proceed by contradiction. Suppose that no such vertex $z$ exists. That is, for every vertex $z\in \{v_2,v_3,\dots,v_{k-1},x,y\}$, we have $z\in N_G(u)$ and $z\in N_G(v_1)$. Let
		\[
		R = N_G(u) \setminus \{v_2, v_3, \dots, v_{k-1}\}, \quad Q = W \setminus N_G(u).
		\]
		By the assumption for contradiction, $R \setminus \{x,y\} \subseteq N_G(v_1)$. Since $x,y$ are common neighbors of $u$ and $v_1$, we have $R \subseteq N_G(v_1)$.
		Moreover, since $G$ is $K_{1,3}$-free, every vertex in $Q$ must be adjacent to $v_1$, which means that $Q \subseteq N_G(v_1)$.
		Thus $N_G(v_1)$ contains both $R$ and $Q$. Consequently,
		\begin{equation}\label{eq:W-upper}
			|W| = |R| + |Q| + (k-2) \le d_G(v_1) + k-2 = k+1 + k-2 = 2k-1.
		\end{equation}
		
		Since $|V(G)| = |W| + |J|$, we have
		\begin{equation}\label{eq:J-lower}
			|J| = |V(G)| - |W| \ge 4k-1 - (2k-1) = 2k.
		\end{equation}
		
		We now count the edges between $W$ and $J$. Since $G$ is $(k+1)$-connected, $e(W,J) \ge (k+1)|J|$. On the other hand, the $K_{1,3}$-free condition implies that $e(W,J) \le 2|W|$. Combining these yields
		\begin{equation}\label{eq:handshake}
			(k+1)|J| \le 2|W|.
		\end{equation}
		
		Substituting \eqref{eq:W-upper} and \eqref{eq:J-lower} into \eqref{eq:handshake}, we obtain
		\[
		2k(k+1) \le 2(2k-1),
		\]
		which simplifies to $k^2 - k + 1 \le 0$. This inequality holds for no positive integer $k$, contradicting the assumption $k \ge 2$. Hence the claim holds.
	\end{claimproof}
	
	By Claim \ref{claim4}, if $d_G(v_1)=k+1$ and $d_G(u)\ge k+2$, there must exist a vertex $z$ such that $zu\in E(G)$ and $zv_1\notin E(G)$, and then we may take $C_k=uxzu$. Set $G'=G-C_k$, then $G'$ is a $k$-connected \(K_{1,3}\)-free split graph.
	
		\begin{submycase}\label{3.2case2.2}
	$u,v_1,v_2,\dots,v_{k-1}\in W$.
\end{submycase}
	
	Since $N_G(u)\neq \emptyset$, let $x\in N_J(u)$. As $d_G(x)\ge k+1$, we may choose $y\in N_G(x)\setminus\{u,v_1,v_2,\dots,\\v_{k-1}\}$. Now consider the common neighbors of $u$ and $y$ in $J$ other than $x$. If there are no other common neighbors, let $C_k=uxyu$. If there exists $j\in N_J(y)\cap N_G(u)$, then let $C_k=uxyju$. Set $G'=G-C_k$, then $G'$ is a $k$-connected \(K_{1,3}\)-free split graph.
	
	By the induction hypothesis, there exist $k-1$ vertex-disjoint cycles $ C_1,C_2,\dots,C_{k-1} $ in $G'$ such that $v_i\in V(C_i)$ for each $i\in [k-1]$ in all the above cases.
	Then we obtain $k$ vertex-disjoint cycles $C_1,C_2,\dots,C_k$ such that $u\in V(C_k)$, $v_i\in V(C_i)$ for each $i\in [k-1]$ and $V(C_1) \cup V(C_2) \cup \cdots \cup V(C_{k}) = V(G)$.
	\end{proof}
	
	\subsection{Edge Spanning Cyclability of $K_{1,3}$-Free Split Graphs }
	
	\begin{mylemma} (Liu et al. \cite{06}) \label{lem3.3}
		 Let $G$ be a $K_{1,3}$-free split graph with a split partition $(W, J)$. Let $a, b \in J$ such that $N_G(a) \cap N_G(b) \neq \emptyset$. Then for any two vertices $x, y \in J \setminus \{a, b\}$, $N_G(x) \cap N_G(y) = \emptyset$.
	\end{mylemma}
	
	\begin{mylemma} \label{lem3.4}
		 Let \(G\) be a \(K_{1,3}\)-free split graph with a split partition \((W,J)\).
		If \(G\) is \(4\)-connected, then \(G\) is \(2\)-edge spanning cyclable.
	\end{mylemma}
	
	\begin{proof}
	Let $G$ be a $4$-connected $K_{1,3}$-free split graph with a split partition $(W,J)$. Take any $e_1,e_2\in E(G)$, and denote $e_1=ab$, $e_2=cd$ with $a,b,c,d\in V(G)$. We shall prove that there exist two vertex-disjoint cycles $C_1$ and $C_2$ such that $e_1\in E(C_1)$, $e_2\in E(C_2)$, and $V(C_1)\cup V(C_2)=V(G)$.
	
	\begin{mycase}\label{3.4case1}
		For any $v\in W$, $|N_J(v)|\le 1$.
	\end{mycase}
	
			\begin{submycase}\label{3.4case1.1}
		$a\in J$ and $b \in W$.
	\end{submycase}
	
	Since $G$ is $4$-connected, we may choose $x \in N_G(a) \setminus \{b, c, d\}$. Let $C_1=abxa$ and set $G'=G-C_1$, then $G'$ is a $3$-connected $K_{1,3}$-free split graph.
	
		\begin{submycase}\label{3.4case1.2}
		$a, b, c, d \in W$.
	\end{submycase}
	
	 When $N_J(a)\cap N_J(b)\neq\emptyset$, let $u\in N_J(a)\cap N_J(b)$ and take $C_1=abua$. When $N_J(a)\cap N_J(b)=\emptyset$, by the maximality of $W$ and the $4$-connectedness of $G$, there must exist a vertex $x\in W\setminus\{a,b,c,d\}$. If $N_J(a)\cap N_J(x)=\emptyset$ and $N_J(b)\cap N_J(x)\neq\emptyset$, we may assume $N_J(b)\cap N_J(x)=\{j\}$, then set $C_1=abjxa$. If $N_J(a)\cap N_J(b)\cap N_J(x)=\emptyset$, take $C_1=abxa$. Set $G'=G-C_1$, then $G'$ is a $3$-connected $K_{1,3}$-free split graph.
	
	\begin{mycase}\label{3.4case2}
		There is some $v\in W$ with $|N_J(v)| \ge 2$.
	\end{mycase}
	
	To avoid an induced $K_{1,3}$, the case $|N_J(v)| > 2$ cannot occur, and hence $|N_J(v)| = 2$.
	At this time, every vertex in $W$ is adjacent to at least one vertex in $J$.

\begin{submycase}\label{3.4case2.1}
	$a\in J$ and $b \in W$.
\end{submycase}
	
	By the 4-connectedness of $G$, we may choose $x\in N_G(a)\setminus\{b,c,d\}$. Apart from $a$, when $b$ and $x$ share no other common neighbors in $J$, the situation is the same as Case \ref{3.4case1.1}. Otherwise, let $j$ be an additional common neighbor of $b$ and $x$ in $J$. For $j\notin\{c,d\}$, take $C_1=abjxa$. If $j\in\{c,d\}$, without loss of generality assume that $j=c$. By the maximality of $W$, $|W|\ge 5$, then there exists a $y\in W\setminus\{b,x,d\}$ such that $ay\in E(G)$ and $cy\notin E(G)$. We may then take $C_1=abya$. Set $G'=G-C_1$, then $G'$ is a $3$-connected $K_{1,3}$-free split graph.

	\begin{submycase}\label{3.4case2.2}
		$a, b, c, d \in W$.
	\end{submycase}
	
	If $|N_J(a)\cap N_J(b)|\le 1$, then the argument is the same as Case \ref{3.4case1.2}. Therefore, we mainly consider the case when $|N_J(a)\cap N_J(b)|=2$.
	Let $\{x,y\}=N_J(a)\cap N_J(b)$. Consider the edge $e_2=cd$. By symmetry, if $|N_J(c)\cap N_J(d)|\le 1$, the situation is the same as in Case \ref{3.4case1.2}. Thus, suppose $|N_J(c)\cap N_J(d)|=2$. By Lemma \ref{lem3.3} and the maximality of $W$, we may assume $\{y,z\}=N_J(c)\cap N_J(d)$. When $|V(G)|\le 9$, we have $|W|=6$. There exist $v_1,v_2\in W\setminus\{a,b,c,d\}$ with $N_J(v_1)=N_J(v_2)=\{x,z\}$. Then we take $C_1=abya$ and $C_2=v_1xv_2zdcv_1$, so that $e_1\in C_1$, $e_2\in C_2$, $V(C_1)\cap V(C_2)=\emptyset$, and $V(C_1)\cup V(C_2)=V(G)$.
	
	For $|V|>9$ and $|W|>6$, to avoid the appearance of $K_{1,3}$, we have either $d_G(y)\ge 5$, or both $d_G(x)\ge 5$ and $d_G(z)\ge 5$. If $d_G(y)\ge 5$, we take $C_1=abxa$; if $d_G(x)\ge 5$ and $d_G(z)\ge 5$, we take $C_1=abya$. Set $G'=G-C_1$, then $G'$ is a $3$-connected $K_{1,3}$-free split graph.
	
	In the above cases, by Theorem \ref{Theorem1.6}, $G'$ is Hamilton-connected, so there exists a Hamiltonian path $P$ between $c$ and $d$ in $G'$. Let $C_2=P+e_2$. Now we obtain two vertex-disjoint cycles $C_1$ and $C_2$ such that $e_1\in E(C_1)$, $e_2\in E(C_2)$, and $V(C_1)\cup V(C_2)=V(G)$.
	\end{proof}

\vskip .2cm
	
	\begin{proof}[Proof of Theorem \ref{Theorem1.16}]
	Now let $G$ be a $2k$-connected $K_{1,3}$-free split graph with a split partition $(W,J)$. We proceed by induction on $k$.
	The base case $k=2$ holds by Lemma \ref{lem3.4}.
	Assume that for $k\ge 3$, every $(2k-2)$-connected $K_{1,3}$-free split graph is $(k-1)$-edge spanning cyclable.
	Now let $e_1,e_2,\dots,e_k$ be any $k$ independent edges of $G$ and set $e_k=ab$ with $a,b\in V(G)$.
	We distinguish two cases: (i) $e_k$ is an edge between $W$ and $J$; (ii) all edges $e_1,e_2,\dots,e_k$ are entirely contained in $W$.
	
	We first consider the case (i). Suppose $a \in J$, $b \in W$.
	Then all the neighbors of $a$ are contained in $W$.
	The set consisting of $b$ and all endpoints of the remaining $k-1$ edges has cardinality at most $2(k-1)+1 = 2k-1$.
	Therefore, we can always pick a vertex
	$x \in N_G(a) \setminus \{b, V(e_1), V(e_2), \dots, V(e_{k-1})\}$ and set $C_k = abxa$.
	Set $G'=G-C_k$, then $G'$ is a $(2k-2)$-connected $K_{1,3}$-free split graph.  By the induction hypothesis, there exist $k-1$ vertex-disjoint cycles $C_1,C_2,\dots,C_{k-1}$ such that $e_i\in E(C_i)$ for each $i\in [k-1]$. Then we obtain $k$ vertex-disjoint cycles $C_1,C_2,\dots,C_k$ such that $e_i\in E(C_i)$ for each $i\in [k]$ and $V(C_1) \cup V(C_2) \cup \cdots \cup V(C_{k}) = V(G)$.
	
	We now address the case (ii). Suppose first that $N_J(a)\cap N_J(b)\neq\emptyset$.
	Select $u\in N_J(a)\cap N_J(b)$ and set $C_k=abua$.
	Now consider $N_J(a)\cap N_J(b)=\emptyset$. By the maximality of $W$ and the $2k$-connectedness of $G$, we can choose a vertex $x\in W\setminus\{a,b\}$ such that $x$ is not incident to any edge $e_i$ with $i\in [k-1]$.
	If $N_J(x)\cap N_J(a)\neq\emptyset$ and $N_J(x)\cap N_J(b)=\emptyset$, we claim $|N_J(x)\cap N_J(a)|=1$.
	Otherwise, there exist two distinct vertices $y,z\in N_J(a)\cap N_J(x)$, which implies that $\{a,y,z,b\}$ induces a $K_{1,3}$, a contradiction.
	Thus we let $N_J(a)\cap N_J(x)=\{y\}$ and take $C_k=abxya$.
	The symmetric scenario with $N_J(x)\cap N_J(a)=\emptyset$ and $N_J(x)\cap N_J(b)\neq\emptyset$ is analogous.
	If $N_J(a)\cap N_J(x)=\{y\}$ and $N_J(b)\cap N_J(x)=\{j\}$, we take $C_k=abjxya$. If $N_J(x)\cap N_J(a)=\emptyset$ and $N_J(x)\cap N_J(b)=\emptyset$, take $C_k=abxa$.
	In all the above cases, set $G'=G-C_k$, then $G'$ is a $(2k-2)$-connected $K_{1,3}$-free split graph.  By the induction hypothesis, there exist $k-1$ vertex-disjoint cycles $C_1,C_2,\dots,C_{k-1}$ such that $e_i\in E(C_i)$ for each $i\in [k-1]$. Then we obtain $k$ vertex-disjoint cycles $C_1,C_2,\dots,C_k$ such that $e_i\in E(C_i)$ for each $i\in [k]$ and $V(C_1) \cup V(C_2) \cup \cdots \cup V(C_{k}) = V(G)$.
	\end{proof}

	\subsection{ Spanning Cyclability of \(\{K_{1,4},K_{1,4}+e\}\)-Free Split Graphs}

\vskip .3cm
	
	\begin{proof}[Proof of Theorem \ref{Theorem1.17}]
		Let G be a $(k+1)$-connected \(\{K_{1,4},K_{1,4}+e\}\)-free split graph with a partition \((W,J)\). The proof is by induction on $k$.

\textbf{Base Case}: Every $2$-connected \(\{K_{1,4},K_{1,4}+e\}\)-free split graph $G$ with  \(|W|\ge |J|\ge 8\) is $1$-spanning cyclable.
Suppose, for contradiction, that $G$ is not. Then by Theorem \ref{Theorem 1.2}, there exists a short cycle $C$ in $G$. Let $C=a_1b_1a_2b_2\cdots a_mb_ma_1$, where $m\ge 2$, and $a_i\in W$, $b_i\in J$ for each $i\in [m]$. We first indicate that $2\le m\le 5$. Indeed, if $m\ge 6$, by the definition of a short cycle, we have $d_G(a_i)=2$, so $b_1a_4,b_1a_6,b_2a_4,b_2a_6\notin E(G)$. Since $W$ is a clique, it follows that $\{a_2,b_1,b_2,a_4,a_6\}\cong K_{1,4}+e$, which is a contradiction. Thus, $2\le m\le 5$. Since $|W|\ge 8$, there exist at least two distinct vertices $v_1,v_2\in W\setminus \{a_1,a_2,a_3,a_4,a_5\} $. By the definition of the short cycle, $b_1,b_2$ are not adjacent to $v_1$ or $v_2$. Then $\{a_2,b_1,b_2,v_1,v_2\}\cong K_{1,4}+e$, which is a contradiction.

\textbf{Hypothesis}: Every $k$-connected, \(\{K_{1,4},K_{1,4}+e\}\)-free split graph $G$ with $|W|\ge |J|\ge 3k+2$ is $(k-1)$-spanning cyclable for $k \ge 2$.

\textbf{Inductive Step}: Let $v_1,v_2,\dots,v_k$ be any $k$ vertices of $G$ and set $u=v_k$.
		
	\begin{claim}\label{claim3.3.4}
		Assume $s \in J$, and there exist $a,b \in N_G(s)$ such that for any $j \in J \setminus \{s\}$, if  $|N_G(j) \cap \{a,b\}| = 2$, then $\max\{d_G(s), d_G(j)\} \ge k+2$.
	\end{claim}
	
	\begin{claimproof}
		We proceed by contradiction. Suppose there exists a vertex $j \in J \setminus \{s\}$ such that $|N_G(j) \cap \{a,b\}| = 2$, but $d_G(s) = d_G(j) = k+1$. Since $a,b \in N_G(s) \cap N_G(j)$, we have $|N_G(s) \cap N_G(j)| \ge 2$. Then
		\[
		|N_G(s) \cup N_G(j)| = d_G(s) + d_G(j) - |N_G(s) \cap N_G(j)| \le (k+1) + (k+1) - 2 = 2k.
		\]
		
		Since $|W| \ge 3k+5$, we have
		\[
		\big|W \setminus \big(N_G(s) \cup N_G(j)\big)\big| = |W| - \big|N_G(s) \cup N_G(j)\big| \ge 3k+5 - 2k = k+5 > 2.
		\]
		Thus there exist two distinct vertices $p,q \in W \setminus \big(N_G(s) \cup N_G(j)\big)$ that are nonadjacent to both $s$ and $j$. Then the subgraph induced by $\{a, s, j, p, q\}$ is isomorphic to $K_{1,4}+e$, a contradiction.
	\end{claimproof}
		
	\begin{claim}\label{claim3.3.6}
		Assume $a,b\in W$, and there exist three distinct vertices $x,y,z\in J$ such that both $a$ and $b$ are adjacent to $x,y,z$. Then at least two of $x,y,z$ have degree at least $k+2$.
	\end{claim}
	
	\begin{claimproof}
		Since $a$ and $b$ are common neighbors of every pair among $x,y,z$, each pair satisfies the premise of Claim \ref{claim3.3.4}. Applying Claim \ref{claim3.3.4} to all three pairs yields
		\[
		\max\{d_G(x), d_G(y)\} \ge k+2, \quad \max\{d_G(x), d_G(z)\} \ge k+2, \quad \max\{d_G(y), d_G(z)\} \ge k+2.
		\]
		Thus at least two of $x,y,z$ have degree at least $k+2$.
	\end{claimproof}
		
		\begin{mycase}\label{3.8case1}
			For any $t\in W$, $|N_J(t)|\le 1$.
		\end{mycase}
		
		By $(k+1)$-connectedness, $e(W,J)\ge (k+1)|J|$. Since each vertex in $W$ has at most one neighbor in $J$, we have $e(W,J)\le |W|$. Thus $|W|\ge (k+1)|J|$. For $k\ge 2$, it follows that $|W|>|J|$.
		
				\begin{submycase}\label{3.8case1.1}
			$u\in J$.
		\end{submycase}
		
		Since $G$ is $(k+1)$-connected, there must exist two vertices $x,y\in N_G(u)\setminus\{v_1,v_2,\dots,v_{k-1}\}$. Let $C_k=uxyu$ and set $G'=G-C_k$, then $G'$ is a $k$-connected $\{K_{1,4},K_{1,4}+e\}$-free split graph.
		
			\begin{submycase}\label{3.8case1.2}
			$u,v_1,v_2,\dots,v_{k-1}\in W$.
		\end{submycase}
		
			If $N_J(u)\neq\emptyset$, then there exists a vertex $z\in N_J(u)$. For $d_G(z)\ge k+1$, we may take $a\in N_G(z)\setminus\{u,v_1,v_2,\dots,v_{k-1}\}$ and set $C_k = uazu$.
		If $N_J(u)=\emptyset$, since $G$ is $(k+1)$-connected and $W$ is a maximal clique, it follows that $|W| \ge k+1+1 = k+2$, then there exist two vertices $x,y\in W\setminus\{u,v_1,v_2,\dots,v_{k-1}\}$. If $N_J(x)\cap N_J(y)=\emptyset$, take $C_k=uxyu$. If $N_J(x)\cap N_J(y)\neq\emptyset$, we may assume $N_J(x)\cap N_J(y)=\{j\}$, then take $C_k=uxjyu$.
	  Set $G'=G-C_k$, then $G'$ is a $k$-connected $\{K_{1,4},K_{1,4}+e\}$-free split graph.
		
		\begin{mycase}\label{3.8case2}
			There is some $t \in W$ with $|N_J(t)|=2$.
		\end{mycase}
		
		Since $G$ is $K_{1,4}$-free, every vertex in $W$ has at most $2$ neighbors in $J$, so $e(W,J)\le 2|W|$. By $(k+1)$-connectedness, $e(W,J)\ge (k+1)|J|$, hence $|W|\ge \frac{k+1}{2}|J|$. For $k\ge 2$, we have $|W|>|J|$.
		
			\begin{submycase}\label{3.8case2.1}
			$u\in J$.
		\end{submycase}

		By the $(k+1)$-connectedness of $G$, we may take two vertices $x,y\in N_G(u)\setminus\{v_1,v_2,\dots,v_{k-1}\}$. Excluding $u$, if $x$ and $y$ have no other common neighbors in $J$, we set $C_k = uxyu$.
		Otherwise, let $j$ be an additional common neighbor of $x$ and $y$ in $J$.
		If $j\neq v_i$ for all $1\le i\le k-1$, take $C_k = uxjyu$.
		Now consider the case $j \in \{v_1, v_2, \dots, v_{k-1}\}$, and without loss of generality let $j = v_1$. Then by Claim \ref{claim3.3.4}, $\max\{d_G(u), d_G(v_1)\} \ge k+2$.
		Therefore, if $d_G(v_1) \ge k+2$, we may take $C_k = uxyu$.
		
		\begin{claim}\label{claim3.3.5}
			If $d_G(u) \ge k+2$ and $d_G(v_1) = k+1$, then there exists a vertex $z \notin \{v_2, \dots, v_{k-1}, x, y\}$ such that $z \in N_G(u)$ and $z \notin N_G(v_1)$.
		\end{claim}
		
		\begin{claimproof}
			We proceed by contradiction. Suppose that no such vertex $z$ exists. That is, for every vertex $z \in \{v_2, v_3, \dots, v_{k-1}, x, y\}$, we have $z \in N_G(u)$ and $z \in N_G(v_1)$. Let
			\[
			R = N_G(u) \setminus \{v_2, v_3, \dots, v_{k-1}\}, \quad Q = W \setminus N_G(u).
			\]
			By the assumption for contradiction, $R \setminus \{x,y\} \subseteq N_G(v_1)$. Since $x,y$ are common neighbors of $u$ and $v_1$, we have $R \subseteq N_G(v_1)$.
			Moreover, since $G$ is $\{K_{1,4}, K_{1,4}+e\}$-free and $u,v_1$ share a common neighbor in $W$, there is at most one vertex in $Q$ nonadjacent to both $u$ and $v_1$, which means $|Q \setminus N_G(v_1)| \le 1$.
			Thus $N_G(v_1)$ contains all vertices in $R$ and all but at most one vertex in $Q$.  Consequently,
			\[
			|W| = |R| + |Q| + (k-2) \le d_G(v_1) + 1 + (k-2) = k+1 + k-1 = 2k.
			\]
			This contradicts $|W| \ge 3k+5$.
		\end{claimproof}
		
		Then by Claim \ref{claim3.3.5}, if $d_G(v_1) = k+1$, there must exist $z$ such that $zu \in E(G)$ and $zv_1 \notin E(G)$, and then we may take $C_k = uxzu$. Set $G'=G-C_k$, then $G'$ is a $k$-connected $\{K_{1,4},K_{1,4}+e\}$-free split graph.
		
			\begin{submycase}\label{3.8case2.2}
			$u,v_1,v_2,\dots,v_{k-1}\in W$.
		\end{submycase}

	If $N_J(u) \neq \emptyset$, let $x \in N_J(u)$. Then we may choose $y \in N_G(x) \setminus \{u, v_1, v_2, \dots, v_{k-1}\}$. If $u$ and $y$ have no other common neighbors in $J$ other than $x$, let $C_k = uxyu$. If there exists $j \in N_J(y) \cap N_J(u)$, then let $C_k = uxyju$.
	When $N_J(u) = \emptyset$, select two vertices $p,q \in W \setminus \{u, v_1, v_2, \dots, v_{k-1}\}$. If $|N_J(p) \cap |N_J(q)| \le 1$, the situation is the same as Case \ref{3.8case1.2}. If $p$ and $q$ possess two common neighbors $c,d$, then by Claim \ref{claim3.3.4},  $\max\{d_G(c), d_G(d)\} \ge k+2$. Accordingly, we may take $C_k = upcqu$ or $C_k = updqu$.
	In the above cases, set $G'=G-C_k$, then $G'$ is a $k$-connected $\{K_{1,4},K_{1,4}+e\}$-free split graph.
		
		\begin{mycase}\label{3.8case3}
			There is some $t \in W$ with $|N_J(t)|\ge 3$.
		\end{mycase}
		
		Since $G$ is $K_{1,4}$-free, every vertex in $W$ has at most $3$ neighbors in $J$, and thus $|N_J(t)| = 3$.
		Moreover, each vertex of $W$ is adjacent to at least one vertex in $J$. Otherwise, the subgraph induced by $t$, its three neighbors in $J$, and a vertex in $W$ with no neighbor in $J$ would be isomorphic to $K_{1,4}$, a contradiction.
		In this case, $e(W,J)\le 3|W|$. By $(k+1)$-connectedness, $e(W,J)\ge (k+1)|J|$, hence $|W|\ge \frac{k+1}{3}|J|$.
		If $|W|=|J|$, then $k+1\le 3$, so $k=2$ and the bipartite subgraph $G[W,J]$ is $3$-regular.
		Take any $z\in W$ and let $N_J(z)=\{j_1,j_2,j_3\}$. Each $j_i$ has exactly $3$ neighbors in $W$, so at most $2$ other neighbors besides $z$. Thus at most $1+2\cdot 3=7$ vertices in $W$ are adjacent to any of $j_1,j_2,j_3$.
		Since $|W|\ge 3k+5=11$, there exist $a,b\in W$ nonadjacent to all of $j_1,j_2,j_3$. Then $\{z,a,b,j_1,j_2\}$ induces a $K_{1,4}+e$ with center $z$, a contradiction.
		Hence $|W|>|J|$.

		\begin{submycase}\label{3.8case3.1}
			$u\in J$.
		\end{submycase}
		
	Since $G$ is $(k+1)$-connected, there exist two vertices $x,y \in N_G(u) \setminus \{v_1, v_2, \dots, v_{k-1}\}$. If $|N_J(x) \cap N_J(y)| \le 2$, the situation is the same as Cases \ref{3.8case1.1} and \ref{3.8case2.1}. Now suppose there exist vertices $j,z$ such that $\{u,j,z\} = N_J(x) \cap N_J(y)$.  Then by Claim \ref{claim3.3.6}, at least two of the three vertices $u,j,z$ have degree at least $k+2$.
	In this case, if $d_G(j) \ge k+2$ and $d_G(z) \ge k+2$, take $C_k = uxyu$. If $d_G(u) \ge k+2$, $d_G(j) \ge k+2$ and $d_G(z) = k+1$, take $C_k = uxzyu$ when $z \neq v_i$ for any $1 \le i \le k-1$. If $z \in \{v_1, v_2, \dots, v_{k-1}\}$, assume without loss of generality that $z = v_1$. Then by Claim \ref{claim3.3.5}, there exists $w \in N_G(u) \setminus \{x,y, v_2, v_3, \dots, v_{k-1}\}$ such that $wv_1 \notin E(G)$, and we may take $C_k = uwxu$.
Set $G'=G-C_k$, then $G'$ is a $k$-connected $\{K_{1,4},K_{1,4}+e\}$-free split graph.
	
			\begin{submycase}\label{3.8case3.2}
			$u,v_1,v_2,\dots,v_{k-1}\in W$.
		\end{submycase}
		
		Take $x\in N_J(u)$. For $d_G(x)\ge k+1$, there exists $y\in N_G(x)\setminus\{u,v_1,v_2,\dots,v_{k-1}\}$.
		If $|N_J(u)\cap N_J(y)| \le 2$, the situation is the same as Cases \ref{3.8case1.2} and \ref{3.8case2.2}. Now suppose that $u$ and $y$ have three distinct common neighbors $x$, $j$ and $z$ in $J$. By claim \ref{claim3.3.6}, at least two of $x,j,z$ satisfy that their degrees are at least $k+2$. If $d_G(j)\ge k+2$ and $d_G(z)\ge k+2$, take $C_k=uxyu$. If $d_G(x) \ge k+2$, $d_G(j) \ge k+2$ and $d_G(z) = k+1$, we take $C_1 = uxyzu$. Set $G'=G-C_k$, then $G'$ is a $k$-connected $\{K_{1,4},K_{1,4}+e\}$-free split graph.
		
In all the cases above, by the induction hypothesis, there exist $k-1$ vertex-disjoint cycles $ C_1,C_2,\dots,C_{k-1} $ in $G'$ such that $v_i\in V(C_i)$ for each $i\in [k-1]$.
Then we obtain $k$ vertex-disjoint cycles $C_1,C_2,\dots,C_k$ such that $v_i\in V(C_i)$ for each $i\in [k]$ and $V(C_1) \cup V(C_2) \cup \cdots \cup V(C_{k}) = V(G)$.
\end{proof}
	
	\subsection{ Edge Spanning Cyclability of \(\{K_{1,4},K_{1,4}+e\}\)-Free Split Graphs}

	\begin{mylemma} \label{lem3.9}
		Let \(G\) be a  $4$-connected \(\{K_{1,4},K_{1,4}+e\}\)-free split graph of order at least $17$, then \(G\) is 2-edge spanning cyclable.
	\end{mylemma}
	
	\begin{proof}
		Let $G$ be a $4$-connected \(\{K_{1,4},K_{1,4}+e\}\)-free split graph with a split partition $(W,J)$. Take any $e_1,e_2\in E(G)$, and denote $e_1=ab$, $e_2=cd$ with $a,b,c,d\in V(G)$. We shall prove that there exist two vertex-disjoint cycles $C_1$ and $C_2$ such that $e_1\in E(C_1)$, $e_2\in E(C_2)$, and $V(C_1)\cup V(C_2)=V(G)$.
		
			\begin{claim}\label{claim3.4.1}
			Assume $s \in J$, and there exist $x,y \in N_G(s)$ such that for any $j \in J \setminus \{s\}$, if $|N_G(j) \cap \{x,y\}| = 2$, then $\max\{d_G(s), d_G(j)\} \ge 5$.
		\end{claim}
		
		\begin{claimproof}
			We proceed by contradiction. Suppose there exists a vertex $j \in J \setminus \{s\}$ such that $|N_G(j) \cap \{x,y\}| = 2$, but $d_G(s) = d_G(j) = 4$.
			Since $x,y \in N_G(s) \cap N_G(j)$, we have $|N_G(s) \cap N_G(j)| \ge 2$. Then
			\[
			|N_G(s) \cup N_G(j)| = d_G(s) + d_G(j) - |N_G(s) \cap N_G(j)| \le 4 + 4 - 2 = 6.
			\]
			
			To exclude an induced $K_{1,4}+e$, there can be at most one vertex in $W$ nonadjacent to both $s$ and $j$. Hence
			\[
			|W| \le |N_G(s) \cup N_G(j)| + 1 = 7.
			\]
			
			Since $G$ is $4$-connected $\{K_{1,4}, K_{1,4}+e\}$-free split graph, then $4|J| \le 2|W|$, which implies
			\[
			|J| \le \frac{1}{2}|W| \le \frac{7}{2},
			\]
			and thus $|J| \le 3$.
			Therefore, $|V(G)| = |W| + |J| \le 7 + 3 = 10$, which contradicts $|V(G)| \ge 17$.
		\end{claimproof}
		
			\begin{claim}\label{claim3.4.2}
			Assume $x,y\in W$, and there exist three distinct vertices $s,t,z\in J$ such that both $x$ and $y$ are adjacent to $s,t,z$. Then at least two of $s,t,z$ have degree at least $5$.
		\end{claim}
		
		\begin{claimproof}
			Since $x$ and $y$ are common neighbors of every pair among $s,t,z$, each pair satisfies the premise of Claim \ref{claim3.4.1}. Applying Claim \ref{claim3.4.1} to all three pairs, we obtain
			\[
			\max\{d_G(s), d_G(t)\} \ge 5, \quad \max\{d_G(s), d_G(z)\} \ge 5, \quad \max\{d_G(t), d_G(z)\} \ge 5.
			\]
			Thus at least two of $s,t,z$ have degree at least $k+2$.
		\end{claimproof}
		
	\begin{mycase}\label{3.9case1}
		For any $v\in W$, $|N_J(v)|\le 1$.
	\end{mycase}
	
	\begin{submycase}\label{3.9case1.1}
		$a\in J$ and $b \in W$.
	\end{submycase}
	
	Since $G$ is $4$-connected, we can choose $x \in N_G(a) \setminus \{b, c, d\}$. Let $C_1=abxa$ and set $G'=G-C_1$, then $G'$ is a $3$-connected \(\{K_{1,4},K_{1,4}+e\}\)-free split graph.
	
		\begin{submycase}\label{3.9case1.2}
		$a, b, c, d \in W$.
	\end{submycase}
	
		When $N_J(a)\cap N_J(b)\neq\emptyset$,
	let $\{j\}=N_J(a)\cap N_J(b)$, we obtain the cycle $C_1=abja$.
	When $N_J(a)\cap N_J(b)=\emptyset$, by the maximality of $W$ and the $4$-connectedness of $G$, there exists $x\in W\setminus\{a,b,c,d\}$.
	If $N_J(a)\cap N_J(x)\cap N_J(b)=\emptyset$, set $C_1=abxa$.
	If $N_J(a)\cap N_J(x)\neq\emptyset$ and $N_J(x)\cap N_J(b)=\emptyset$, let $\{y\}=N_J(a)\cap N_J(x)$.
	We then take $C_1=abxya$.
Set $G'=G-C_1$, then $G'$ is a $3$-connected \(\{K_{1,4},K_{1,4}+e\}\)-free split graph.
	
	\begin{mycase}\label{3.9case2}
		There is some $v \in W$ with $|N_J(v)|=2$.
	\end{mycase}
	
	In this case, at most one vertex in $W$ has no neighbor in $J$. Otherwise, an induced $K_{1,4}+e$ appears, a contradiction.
	
			\begin{submycase}\label{3.9case2.1}
		$a\in J$ and $b \in W$.
	\end{submycase}

As $G$ is $4$-connected, we may take $x \in N_G(a) \setminus \{b,c,d\}$.
 The case $|N_J(x) \cap N_J(b)| = 1$ is analogous to Case \ref{3.9case1.1}. We next consider the case when $|N_J(x) \cap N_J(b)| = 2$. Let $j \in N_J(x) \cap N_J(b) \setminus \{a\}$. If $j \neq c$, we take $C_1 = abjxa$. If $j = c$, then by Claim \ref{claim3.4.1}, $\max\{d_G(a), d_G(c)\} \ge 5$. Now if $d_G(c) \ge 5$, we set $C_1 = abxa$. If $d_G(a) \ge 5$ and $d_G(c) = 4$, there exists $y \in W$ satisfying $cy \notin E(G)$ and $ay \in E(G)$. Accordingly, we may take $C_1 = abya$. Set $G'=G-C_1$, then $G'$ is a $3$-connected \(\{K_{1,4},K_{1,4}+e\}\)-free split graph.

		\begin{submycase}\label{3.9case2.2}
		$a, b, c, d \in W$.
	\end{submycase}
	
When $N_J(a) \cap N_J(b) \neq \emptyset$, the situation $|N_J(a) \cap N_J(b)| \le 1$ is the same as Case \ref{3.9case1.2}.  Now if $|N_J(a) \cap N_J(b)| = 2$, let $\{x,y\} = N_J(a) \cap N_J(b)$. By Claim \ref{claim3.4.1}, we obtain $d_G(x) \ge 5$ or $d_G(y) \ge 5$. We then take $C_1 = abya$ or $C_1 = abxa$, respectively.
When $N_J(a) \cap N_J(b) = \emptyset$, by the maximality of $W$ and the 4-connectedness of $G$, there exists a vertex $z \in W \setminus \{a,b,c,d\}$. The subcases $N_J(a) \cap N_J(b) \cap N_J(z) = \emptyset$, or $|N_J(a) \cap N_J(z)| = 1$ together with $N_J(z) \cap N_J(b) = \emptyset$, are analogous to Case \ref{3.9case1.2}.
We next treat $|N_J(a) \cap N_J(z)| = 2$ under the assumption $N_J(z) \cap N_J(b) = \emptyset$. Let $\{p,q\} = N_J(a) \cap N_J(z)$. By Claim \ref{claim3.4.1}, at least one of $d_G(p) \ge 5$ or $d_G(q) \ge 5$ holds. If $d_G(p) \ge 5$, set $C_1 = abzqa$, if $d_G(q) \ge 5$, set $C_1 = abzpa$. The symmetric argument applies for $|N_J(b) \cap N_J(z)| = 2$ and $N_J(z) \cap N_J(a) = \emptyset$.

Now suppose $|N_J(a)\cap N_J(z)| = 1$ and $|N_J(b)\cap N_J(z)| = 1$.
Let $\{g\} = N_J(a)\cap N_J(z)$ and $\{h\} = N_J(b)\cap N_J(z)$, so $|N_G(g)\cap N_G(h)|\ge 1$.
To avoid an induced $K_{1,4}+e$, at most one vertex in $W$ is nonadjacent to both $g$ and $h$, which gives
\[
|N_G(g)\cup N_G(h)| \ge |W| - 1.
\]
At this time, we claim that $\max\{d_G(g), d_G(h)\} \ge 5$.
Suppose for contradiction that $d_G(g) = d_G(h) = 4$. Then
\[
|N_G(g)\cup N_G(h)| = d_G(g)+d_G(h) - |N_G(g)\cap N_G(h)| \le 4+4 - 1 = 7.
\]
Combining the two inequalities above yields $|W| - 1 \le 7$, so $|W| \le 8$.
Then $|J| \le \frac{1}{2}|W| \le 4$, and hence $|V(G)| = |W|+|J| \le 12$, which contradicts $|V(G)|\ge 17$.
Therefore $\max\{d_G(g), d_G(h)\} \ge 5$.
Now if $d_G(g)\ge 5$, set $C_1 = abzha$, if $d_G(h)\ge 5$, set $C_1 = abzga$.
In the above cases, set $G'=G-C_1$, then $G'$ is a $3$-connected \(\{K_{1,4},K_{1,4}+e\}\)-free split graph.	
	\begin{mycase}\label{3.9case3}
		There is some $v \in W$ with $|N_J(v)| \ge 3$.
	\end{mycase}
	To avoid an induced \(\{K_{1,4},K_{1,4}+e\}\), the case $|N_J(v)| > 3$ cannot occur, and hence $|N_J(v)| = 3$.
	In this case, every vertex in $W$ is adjacent to at least one vertex in $J$.
	
	\begin{submycase}\label{3.9case3.1}
	$a\in J$ and $b \in W$.
\end{submycase}
	
	It follows from 4-connectedness that we can choose a vertex $x \in N_G(a) \setminus \{b,c,d\}$. The case $|N_J(x) \cap N_J(b)| \le 2$ is analogous to Cases \ref{3.9case1.1} and \ref{3.9case2.1}.
	We next consider $|N_J(x) \cap N_J(b)| = 3$. Let $y,z \in N_J(x) \cap N_J(b) \setminus \{a\}$. Then by Claim \ref{claim3.4.2}, at least two among $a,y,z$ have degree at least 5.	Now if $d_G(y) \ge 5$, $d_G(z) \ge 5$, set $C_1 = abxa$. If $d_G(a) \ge 5$, $d_G(y) \ge 5$ and $d_G(z) = 4$, we take $C_1 = abzxa$ when $z \neq c$. If $z = c$, there must exist a vertex $w$ satisfying $aw \in E(G)$ and $cw \notin E(G)$. Then we take $C_1=abwa$. Set $G'=G-C_1$, then $G'$ is a $3$-connected \(\{K_{1,4},K_{1,4}+e\}\)-free split graph.

		\begin{submycase}\label{3.9case3.2}
	$a, b, c, d \in W$.
\end{submycase}

When $N_J(a)\cap N_J(b) \neq \emptyset$, the subcase $|N_J(a)\cap N_J(b)| \le 2$ is analogous to Cases \ref{3.9case1.2} and \ref{3.9case2.2}. For $|N_J(a)\cap N_J(b)| = 3$, let $\{s,t,u\} = N_J(a)\cap N_J(b)$. By Claim \ref{claim3.4.2}, at least two of $s,t,u$ have degree at least 5. Without loss of generality, let $d_G(t) \ge 5$ and $d_G(u) \ge 5$. We set $C_1 = absa$.
When $N_J(a)\cap N_J(b) = \emptyset$, the maximality of $W$ and the $4$-connectedness of $G$ yield a vertex $x\in W\setminus\{a,b,c,d\}$. The subcases $|N_J(a)\cap N_J(x)| \le 1$ with $|N_J(b)\cap N_J(x)| \le 1$ and $|N_J(a)\cap N_J(x)| = 2$ with $|N_J(b)\cap N_J(x)| =\emptyset$ are analogous to Cases \ref{3.9case1.2} and \ref{3.9case2.2}.
Note that $|N_J(a)\cap N_J(x)|$ cannot be 3. Otherwise, let $\{p,q,r\} = N_J(a)\cap N_J(x)$. Since $N_J(a)\cap N_J(b) = \emptyset$, $b$ is nonadjacent to all of $p,q,r$. Then the subgraph induced by $\{x,p,q,r,b\}$ is isomorphic to $K_{1,4}$, a contradiction.
Now suppose that $|N_J(a) \cap N_J(x)| = 2$ and $|N_J(x) \cap N_J(b)| = 1$, let $\{u,f\} = N_J(a) \cap N_J(x)$ and $\{j\} = N_J(b) \cap N_J(x)$. Since $N_J(a) \cap N_J(b) = \emptyset$, both $u$ and $f$ are nonadjacent to $b$. Then the subgraph induced by $\{x, u, f, j, b\}$ is isomorphic to $K_{1,4}+e$, a contradiction, thus this case cannot occur.
Set $G'=G-C_1$, then $G'$ is a $3$-connected \(\{K_{1,4},K_{1,4}+e\}\)-free split graph.

In the above cases, by Theorem \ref{Theorem1.12}, $G'$ is Hamilton-connected, so there exists a Hamiltonian path $P$ between $c$ and $d$ in $G'$. Let $C_2=P+e_2$.
Now we obtain two vertex-disjoint cycles $C_1$ and $C_2$ such that $e_1\in E(C_1)$, $e_2\in E(C_2)$, and $V(C_1)\cup V(C_2)=V(G)$.
\end{proof}

\vskip .3cm	

	\begin{proof}[Proof of Theorem \ref{Theorem1.19}]
		Now let $G$ be a $2k$-connected $\{K_{1,4},K_{1,4}+e\}$-free split graph with a split partition $(W,J)$.
		We proceed by induction on $k$.
		The base case $k=2$ follows from Lemma \ref{lem3.9}.
		Assume that for $k\ge 3$, every $(2k-2)$-connected $\{K_{1,4},K_{1,4}+e\}$-free split graph is $(k-1)$-edge spanning cyclable.
		Now let $e_1,e_2,\dots,e_k$ be any $k$ independent edges of $G$ and set $e_k=ab$ with $a,b\in V(G)$.
		We distinguish two cases: (i) $e_k$ is an edge between $W$ and $J$; (ii) all edges $e_1,e_2,\dots,e_k$ are entirely contained in $W$.
		
		We first treat the case (i). Without loss of generality, let $a\in J$ and $b\in W$.
		As $a$ belongs to $J$, all its neighbors are contained in $W$.
		The set consisting of $b$ and all endpoints of the remaining $k-1$ edges has cardinality at most $2(k-1)+1 = 2k-1$.
		Therefore, we can always pick a vertex
		$x \in N_G(a) \setminus \{b, V(e_1), V(e_2), \dots, V(e_{k-1})\}$ and set $C_k = abxa$.
		Set $G'=G-C_k$, then $G'$ is a $(2k-2)$-connected $\{K_{1,4},K_{1,4}+e\}$-free split graph.  By the induction hypothesis, there exist $k-1$ vertex-disjoint cycles $C_1,C_2,\dots,C_{k-1}$ such that $e_i\in E(C_i)$ for each $i\in [k-1]$. Then we obtain $k$ vertex-disjoint cycles $C_1,C_2,\dots,C_k$ such that $e_i\in E(C_i)$ for each $i\in [k]$ and $V(C_1) \cup V(C_2) \cup \cdots \cup V(C_{k}) = V(G)$.
		
		We now turn to the case (ii).
		Suppose first that $N_J(a)\cap N_J(b)\neq\emptyset$.
		Select $u\in N_J(a)\cap N_J(b)$ and set $C_k=abua$.
		Now suppose that $N_J(a)\cap N_J(b)=\emptyset$. By the maximality of $W$ and the $2k$-connectedness of $G$, we can choose a vertex $x\in W\setminus\{a,b\}$ such that $x$ is not incident to any edge $e_i$ with $i\in [k-1]$.
	If $N_J(a) \cap N_J(x) = \emptyset$ and $N_J(b) \cap N_J(x) = \emptyset$, we take $C_k = abxa$. If $N_J(a) \cap N_J(x) \neq \emptyset$ while $N_J(b) \cap N_J(x) = \emptyset$, let $j \in N_J(a) \cap N_J(x)$ and take $C_k = abxja$. By symmetry, if $N_J(a) \cap N_J(x) = \emptyset$ and $N_J(b) \cap N_J(x) \neq \emptyset$, let $y \in N_J(b) \cap N_J(x)$ and take $C_k = abyxa$. If $N_J(a) \cap N_J(x)\neq \emptyset $ and $N_J(b) \cap N_J(x)\neq \emptyset$, let $p \in N_J(a) \cap N_J(x)$ and $q \in N_J(b) \cap N_J(x)$. Then we take $C_k=abxpa$ or $C_k=abqxa$.
	In the above cases, set $G'=G-C_k$, then $G'$ is a $(2k-2)$-connected $\{K_{1,4},K_{1,4}+e\}$-free split graph.  By the induction hypothesis, there exist $k-1$ vertex-disjoint cycles $C_1,C_2,\dots,C_{k-1}$ such that $e_i\in E(C_i)$ for each $i\in [k-1]$. Then we obtain $k$ vertex-disjoint cycles $C_1,C_2,\dots,C_k$ such that $e_i\in E(C_i)$ for each $i\in [k]$ and $V(C_1) \cup V(C_2) \cup \cdots \cup V(C_{k}) = V(G)$.
	\end{proof}

\end{document}